\documentclass[11pt]{article}

\usepackage{amscd,amsmath, amssymb, fancyhdr, mathbbol}

\numberwithin{equation}{section}

\newcommand{\version}{version 3.0,\ \   July 06, 2014}

\def\eqref#1{(\ref{#1})}

\newcommand{\arrow}{{\:\longrightarrow\:}}
\newcommand{\Z}{{\mathbb Z}}
\def\C{{\mathbb C}}
\newcommand{\R}{{\mathbb R}}
\newcommand{\Q}{{\mathbb Q}}

\def\1{\sqrt{-1}\:}
\newcommand{\restrict}[1]{{\left|_{{\phantom{|}\!\!}_{#1}}\right.}}
\newcommand{\cntrct}                
{\hspace{2pt}\raisebox{1pt}{\text{$\lrcorner$}}\hspace{2pt}}

\renewcommand{\tilde}{\widetilde}
\renewcommand{\bar}{\overline}
\renewcommand{\phi}{\varphi}
\renewcommand{\epsilon}{\varepsilon}
\renewcommand{\geq}{\geqslant}
\renewcommand{\leq}{\leqslant}

\newcommand{\Gr}{\operatorname{Gr}}
\newcommand{\Teich}{\operatorname{Teich}}
\newcommand{\Comp}{\operatorname{\sf Comp}}
\newcommand{\Per}{\operatorname{\sf Per}}
\newcommand{\Perspace}{\operatorname{{\mathbb P}\sf er}}

\newcommand{\im}{\operatorname{im}}
\newcommand{\End}{\operatorname{End}}

\newcommand{\Id}{\operatorname{Id}}

\newcommand{\Sym}{\operatorname{Sym}}

\newcommand{\Diff}{\operatorname{Diff}}

\newcommand{\rk}{\operatorname{rk}}

\newcommand{\Tw}{\operatorname{Tw}}

\renewcommand{\Re}{\operatorname{Re}}
\renewcommand{\Im}{\operatorname{Im}}

\newcounter{Mycounter}[section]
\newcounter{lemma}[section]
\renewcommand{\thelemma}{{Lemma \thesection.\arabic{lemma}}}
\newcommand{\lemma}{%
    \setcounter{lemma}{\value{Mycounter}}
    \refstepcounter{lemma}
    \stepcounter{Mycounter}
    {\noindent \bf \thelemma:\ }}

\newcounter{claim}[section]
\renewcommand{\theclaim}{{Claim \thesection.\arabic{claim}}}
\newcommand{\claim}{%
    \setcounter{claim}{\value{Mycounter}}
    \refstepcounter{claim}
    \stepcounter{Mycounter}
    {\noindent \bf \theclaim:\ }}

\newcounter{sublemma}[section]
\renewcommand{\thesublemma}{{Sublemma \thesection.\arabic{sublemma}}}
\newcommand{\sublemma}{%
    \setcounter{sublemma}{\value{Mycounter}}
    \refstepcounter{sublemma}
    \stepcounter{Mycounter}
    {\noindent \bf \thesublemma:\ }}

\newcounter{corollary}[section]
\renewcommand{\thecorollary}{{Corollary \thesection.\arabic{corollary}}}
\newcommand{\corollary}{%
    \setcounter{corollary}{\value{Mycounter}}
    \refstepcounter{corollary}
    \stepcounter{Mycounter}
    {\noindent \bf \thecorollary:\ }}

\newcounter{theorem}[section]
\renewcommand{\thetheorem}{{Theorem \thesection.\arabic{theorem}}}
\newcommand{\theorem}{%
    \setcounter{theorem}{\value{Mycounter}}
    \refstepcounter{theorem}
    \stepcounter{Mycounter}
    {\noindent \bf \thetheorem:\ }}

\newcounter{conjecture}[section]
\renewcommand{\theconjecture}{{Conjecture \thesection.\arabic{conjecture}}}
\newcommand{\conjecture}{%
    \setcounter{conjecture}{\value{Mycounter}}
    \refstepcounter{conjecture}
    \stepcounter{Mycounter}
    {\noindent \bf \theconjecture:\ }}

\newcounter{proposition}[section]
\renewcommand{\theproposition}
      {{Proposition \thesection.\arabic{proposition}}}
\newcommand{\proposition}{%
    \setcounter{proposition}{\value{Mycounter}}
    \refstepcounter{proposition}
    \stepcounter{Mycounter}
    {\noindent \bf \theproposition:\ }}

\newcounter{definition}[section]
\renewcommand{\thedefinition}
      {{Definition~\thesection.\arabic{definition}}}
\newcommand{\definition}{%
    \setcounter{definition}{\value{Mycounter}}
    \refstepcounter{definition}
    \stepcounter{Mycounter}
    {\noindent \bf \thedefinition:\ }}

\newcounter{example}[section]
\newcounter{remark}[section]
\renewcommand{\theremark}{{Remark \thesection.\arabic{remark}}}
\newcommand{\remark}{%
    \setcounter{remark}{\value{Mycounter}}
    \refstepcounter{remark}
    \stepcounter{Mycounter}
    {\noindent \bf \theremark:\ }}

\newcounter{problem}[section]
\newcounter{question}[section]
\makeatletter

\@addtoreset{equation}{section} \@addtoreset{footnote}{section}
\makeatother

\def\blacksquare{\hbox{\vrule width 5pt height 5pt depth 0pt}}
\def\endproof{\blacksquare}

\begin{document}
\begin{center}
{\LARGE\bf
Degenerate twistor spaces for hyperk\"ahler manifolds\\[4mm]
}

 Misha
Verbitsky\footnote{Partially supported by RFBR grants
 12-01-00944-Á,  10-01-93113-NCNIL-a, and
AG Laboratory NRI-HSE, RF government grant, ag. 11.G34.31.0023, and
the Simons-IUM fellowship grant.}

\end{center}

{\small \hspace{0.1\linewidth}
\begin{minipage}[t]{0.8\linewidth}
{\bf Abstract} \\
Let $M$ be a hyperk\"ahler manifold, and $\eta$ 
a closed, positive (1,1)-form with $\rk \eta < \dim M$.
We associate to $\eta$ a family of complex structures on $M$,
called a degenerate twistor family, and parametrized by 
a complex line. When $\eta$ is a pullback of a K\"ahler form
under a Lagrangian fibration $L$, all the fibers of degenerate twistor
family also admit a Lagrangian fibration, with the fibers isomorphic
to that of $L$. Degenerate twistor families can be obtained by taking limits of
twistor families, as one of the K\"ahler forms in the hyperk\"ahler
triple goes to $\eta$.
\end{minipage}
}

{\scriptsize
\tableofcontents
}


\section{Introduction}


\subsection[Complex structures obtained from non-degenerate closed 2-forms]{Complex structures obtained from non-degenerate\\ closed 2-forms}

The degenerate twistor spaces (\ref{_dege_twi_Definition_})
are obtained through the following construction.

\hfill

\definition\label{_non-dege_form_Definition_}
A complex-valued 2-form $\Omega$ on a real manifold
$M$ is called {\bf non-degenerate} if $\Omega(v, \cdot)\neq 0$
for any non-zero tangent vector $v\in T_m M$.
Complex structures on $M$ can be obtained from complex
sub-bundles $B=T^{1,0}M\subset TM\otimes_\R \C$
satisfying 
\begin{equation}\label{_complex_str_via_sub-bu_Equation_}
B\oplus \bar B=TM\otimes_\R \C, \ \ [B,B]\subset B
\end{equation}
(\ref{_almost_co_by_(0,1)_Claim_}).

To obtain such $B$, take a non-degenerate 
(\ref{_non-dege_form_Definition_}), 
closed 2-form
$\Omega\in \Lambda^2(M,\C)$, satisfying 
$\Omega^{n+1}=0$, where $4n=\dim_\R M$. Then
$\ker\Omega:=\{v\in T_m M\otimes_\R \C \ \ |\ \ \Omega(v, \cdot)=0\}$
satisfies the conditions of \eqref{_complex_str_via_sub-bu_Equation_}
(see \ref{_comple_stru_from_2-form_Theorem_}).

Degenerate twistor spaces are obtained by constructing
a family $\Omega_t$ of such 2-forms, parametrized by $t\in \C$,
on hyperk\"ahler manifolds. The relation $\Omega_t^{n+1}=0$ 
 follows from the properties of cohomology
of hyperk\"ahler manifolds, most notably the Fujiki formula, 
computation of cohomology performed in \cite{_Verbitsky:cohomo_}, 
and positivity (see Subsection \ref{_dege_twi_Subsection_}).

\subsection{Degenerate twistor families and Teichm\"uller spaces}
\label{_dege_twi_Teich_Subsection_}

In this subsection, we provide a motivation for
the term ``degenerate twistor family''. We introduce the
twistor families of complex structures on 
hyperk\"ahler manifolds and the corresponding
rational curves in the moduli, called {\bf the twistor lines}.

A degenerate twistor family is a family ${\cal Z}$
of deformations of a holomorphically symplectic
manifold $(M,\Omega)$ associated with a positive, closed, semidefinite
form $\eta$ satisfying $\eta^{n-i}\wedge \Omega^{i+1}=0$, for all
$i=0,1, ..., n$, where
$\dim_\C M=2n$ (\ref{_semipis_integra_Theorem_}). 
In this subsection, we define a twistor family of
a hyperk\"ahler manifold, and explain 
how these families can be obtained as limits of twistor deformations.

Throughout this paper, {\bf a hyperk\"ahler manifold}
is a compact, holomorphically symplectic manifold $M$ of 
K\"ahler type. It is called {\bf simple} (\ref{_hk_simple_Definition_})
if $\pi_1(M)=0$ and $H^{2,0}(M)=\C$. We shall (sometimes silently) 
assume that all hyperk\"ahler manifolds we work with are simple.

A hyperk\"ahler metric is a metric $g$ compatible with three
complex structures $I,J,K$ satisfying the quaternionic relations
$IJ=-JI=K$, which is K\"ahler with respect to $I,J,K$.
By the Calabi-Yau theorem, any compact, holomorphically symplectic manifold of 
K\"ahler type admits a hyperk\"ahler metric, which is unique
in each K\"ahler class (\ref{_Calabi-Yau_Theorem_}).

{\bf A hyperk\"ahler structure} is a hyperk\"ahler metric $g$
together with the compatible quaternionic action, that is,
a triple of complex structures satisfying the quaternionic
relations and K\"ahler.  For any $(a, b, c)\in S^2\subset \R^3$,
the quaternion $L:=aI+bJ+cK$ defines another complex
structure on $M$, also K\"ahler with respect to $g$.
This can be seen because the Levi-Civita connection $\nabla$
of $(M,g)$ preserves $I,J,K$, hence $\nabla L=0$,
and this implies integrability and K\"ahlerness of $L$.

Such a complex structure is called {\bf induced complex structure}.
The $\C P^1$-family of induced complex structures obtained this way 
is in fact holomorphic (Subsection \ref{_hk_basic_Subsection_}).
It is called {\bf the twistor deformation}. 
The twistor families can be described in terms of
periods of hyperk\"ahler manifolds as follows.

\hfill

\definition
Let $M$ be a compact complex manifold, and 
$\Diff_0(M)$ a connected component of its diffeomorphism group
(also known as {\bf the group of isotopies}). Denote by $\Comp$
the space of complex structures on $M$, equipped with
topology induced from the $C^\infty$-topology on the space
of all tensors, and let $\Teich:=\Comp/\Diff_0(M)$. We call 
it {\bf the Teichm\"uller space.}

\hfill

\definition
 Let 
\[ \Per:\; \Teich \arrow {\mathbb P}H^2(M, \C)
\]
map $J$ to a line $H^{2,0}(M,J)\in {\mathbb P}H^2(M, \C)$.
The map $\Per$ is called {\bf the period map}.

\hfill

For a simple hyperk\"ahler manifold, an important
bilinear symmetric form $q\in \Sym^2H^2(M,\Q)^*$
is defined, called {\bf Bogomolov-Beauville-Fujiki form}
(\ref{_BBF_Definition_}). This form is a topological invariant
of the manifold $M$, allowing one to describe deformations
of a complex structure very explicitly. Recall that two
points $x, y$ on a topological space are called {\bf non-separable},
if all their neighbourhoods $U_x\ni x$, $U_y\ni y$ intersect.
We denote the corresponding symmetric relation in
$\Teich$ by $x\sim y$. D. Huybrechts has shown that
$x\sim y$ for $x, y\in \Teich$ 
implies that the corresponding complex 
manifolds $(M,x)$ and $(M,y)$ are bimeromorphic
(\cite{_Huybrechts:basic_}). In \cite{_V:Torelli_}
it was shown that $\sim$ defines an equivalence
relation on $\Teich$; the corresponding quotient
space $\Teich/\sim$ is called {\bf the birational
Teichm\"uller space}, and denoted $\Teich_b$. 

Define the {\bf period space} $\Perspace$ as 
\[ \Perspace:=\{l\in {\mathbb P}(H^2(M,\C))\ \ |\ \ q(l,l)=0, q(l, \bar l)>0\}.
\]
The global Torelli theorem (\cite{_V:Torelli_})
can be stated as follows.

\hfill

\theorem\label{_Global_Torelli_Theorem_}
Let $M$ be a simple hyperk\"ahler manifold, $\Teich_b$
the birational Teichm\"uller space, and
$\Per:\; \Teich_b \arrow  {\mathbb P}(H^2(M,\C))$
the period map. Then $\Per$ maps $\Teich_b$ to $\Perspace$,
inducing a diffeomorphism of each connected component 
of $\Teich_b$ with $\Perspace$.

{\bf Proof:} See \cite{_V:Torelli_}. \endproof

\hfill

\remark\label{_Perspace_Grasmann_Remark_}
The period space $\Perspace$ is equipped with a transitive action
of $SO(H^2(M,\R))$. Using this action, one can identify
$\Perspace$ with the Grassmann space  of 2-dimensional,
positive, oriented planes
$\Gr_{{}_{+,+}}(H^2(M,\R))=SO(b_2-3,3)/SO(2) \times SO(b_2-3,1)$.
Indeed, for each $l\in {\mathbb P}H^2(M, \C)$, the space
generated by $\langle \Im l, \Re l\rangle$ is 2-dimensional, because 
$q(l,l)=0, q(l, \bar l)\neq 0$ implies that $l \cap H^2(M,\R)=0$.
This produces a point of $\Gr_{{}_{+,+}}(H^2(M,\R))$ from
$l\in \Perspace$. To obtain the converse correspondence, notice that
for any 2-dimensional positive
plane  $V\in H^2(M,\R)$, 
the quadric $\{l\in V \otimes_\R \C\ \ |\ \ q(l,l)=0\}$
consists of two lines $l\in \Perspace$. A choice of one of two lines is 
determined by the orientation in $V$.

\hfill

We shall describe the Teichm\"uller space and the moduli
of hyperk\"ahler structures in the same spirit, as follows.

Recall that any hyperk\"ahler structure $(M,I,J,K,g)$
defines a triple of K\"ahler forms 
$\omega_I, \omega_J, \omega_K\in \Lambda^2(M)$
(Subsection \ref{_hk_basic_Subsection_}). A hyperk\"ahler structure
on a simple hyperk\"ahler manifold is determined by a 
complex structure and a K\"ahler class (\ref{_Calabi-Yau_Theorem_}).

We call hyperk\"ahler structures {\bf equivalent}
if they can be obtained by a homothety and a 
quaternionic reparametrization:
\[ (M,I,J,K,g)\sim (M,hIh^{-1},hJh^{-1},hKh^{-1},\lambda g),\]
for $h\in {\mathbb H}^*$, $\lambda\in \R^{>0}$. 
Let $\Teich^{\cal H}$ be the set of equivalence
classes of hyperk\"ahler structrues up to the action
of $\Diff_0(M)$, and $\Teich^{\cal H}_b$ its quotient
by $\sim$ (the non-separability relation).

\hfill

\theorem\label{_hk_Torelli_Theorem_}
Consider the period map 
\[ \Per_{\cal H}:\; \Teich^{\cal H}_b\arrow \Gr_{+++}(H^2(M,\R))\]
associating the plane $\langle \omega_I, \omega_J, \omega_K\rangle$
in the Grassmannian of 3-dimensional positive oriented planes
to an equivalence class of hyperk\"ahler structures.
Then $\Per_{\cal H}$ is injective, and defines
an open embedding on each connected component of
$\Teich^{\cal H}_b$.

\hfill

{\bf Proof:} 
As follows from global Torelli theorem (\ref{_Global_Torelli_Theorem_})
and \ref{_Perspace_Grasmann_Remark_}, a 
complex structure is determined (up to diffeomorphism
and a birational equivalence) by a 2-plane
$V\in \Gr_{{}_{+,+}}(H^2(M,\R))=SO(b_2-3,3)/SO(2) \times SO(b_2-3,1)$,
where $V=\langle \Re\Omega,\Im\Omega\rangle$, and $\Omega$
a holomorphically symplectic form (defined uniquely up to a
multiplier).
Let $\omega\in H^{1,1}(M,I)=V^\bot$ be a K\"ahler form. The corresponding
hyperk\"ahler structure gives an orthogonal triple of K\"ahler forms
$\omega_J, \omega_K\in V, \omega_I:=\omega\in V^\bot$
satisfying $q(\omega_I,\omega_I)=q(\omega_J,\omega_J)=q(\omega_K,\omega_K)=C$.
The group $SU(2)\times \R^{>0}$ acts on the set of such orthogonal
bases transitively. Therefore, a hyperk\"ahler structure
is determined (up to equivalence of hyperk\"ahler structures
and non-separability) by a 3-plane 
$W=\langle \omega_I,\omega_J, \omega_K\rangle\subset H^2(M,\R)$.

We have shown that $\Per_{\cal H}$ is injective.
To finish the proof of \ref{_hk_Torelli_Theorem_},
it remains to show that  $\Per_{\cal H}$ is an open embedding.
However, for a sufficiently small 
$v\in \langle \omega_J, \omega_K\rangle^\bot=H^{1,1}_\R(M,I)$, the 
form $v+\omega_I$ is also K\"ahler (the K\"ahler cone is open
in $H^{1,1}_\R(M,I)$), hence 
$W'=\langle \omega_I+v,\omega_J, \omega_K\rangle$ also belongs to
an image of $\Per_{\cal H}$. This implies that the
differential $D(\Per_{\cal H})$ is surjective.
\endproof

\hfill

Every hyperk\"ahler structure induces a whole 2-dimensional
sphere of complex structures on $M$, as follows. 
Consider a triple $a, b, c\in\R$, $a^2 + b^2+ c^2=1$,
and let $L:= aI + bJ +cK$ be the corresponding quaternion. 
Quaternionic relations imply immediately that $L^2=-1$,
hence $L$ is an almost complex structure. 
Since $I, J, K$ are K\"ahler, they are parallel with respect
to the Levi-Civita connection. Therefore, $L$ is also parallel.
Any parallel complex structure is integrable, and K\"ahler.
We call such a complex structure $L= aI + bJ +cK$
a {\bf complex structure induced by the hyperk\"ahler structure}.
The corresponding complex manifold is denoted by $(M,L)$.
There is a holomorphic family of 
induced complex structures, parametrized by $S^2=\C P^1$.
The total space
of this family is called the \emph{\bf twistor space}
of a hyperk\"ahler manifold; it is constructed as follows.

Let $M$ be a hyperk\"ahler manifold. Consider the product  $\Tw(M) = M
\times S^2$. Embed the sphere 
$S^2 \subset {\mathbb H}$ into the quaternion algebra
${\mathbb H}$ as the set of all quaternions $J$ with $J^2 = -1$. For every point
$x = m \times J \in X = M \times S^2$ the tangent space $T_x\Tw(M)$ is
canonically decomposed $T_xX = T_mM \oplus T_JS^2$. Identify $S^2$
with $\C P^1$, and let $I_J:T_JS^2 \to T_JS^2$ be the complex 
structure operator. Consider the complex structure $I_m:T_mM \to T_mM$ 
on $M$ induced by $J \in S^2 \subset {\mathbb H}$.

The operator $I_{\Tw} = I_m \oplus I_J:T_x\Tw(M) \to T_x\Tw(M)$ 
satisfies $I_{\Tw} \circ I_{\Tw} =
-1$. It depends smoothly on the point $x$, hence it defines 
an almost complex structure on $\Tw(M)$. This almost 
complex structure is known to be integrable
(see e.g. \cite{_Salamon_}, \cite{_Kaledin:twistor_}). 

\hfill

\definition
The space $\Tw(M)$ constructed above is called {\bf the twistor space}
of a hyperk\"ahler manifold.

\hfill

The twistor space defines a family of deformations of a complex
structire on $M$, called {\bf the twistor family}; the corresponding
curve in the Teichm\"uller space is called {\bf the twistor line}.

Let $(M,I,J,K)$ be a hyperk\"ahler structure, and 
$W=\langle \omega_I, \omega_J, \omega_K\rangle$
the corresponding 3-dimensional plane. The twistor
family gives a rational line $\C P^1\subset \Teich$,
which can be recovered from $W$ as follows. Recall that
by global Torelli theorem, each component of $\Teich$
is identified (up to gluing together non-separable points)
with the Grassmannian $\Gr_{{}_{+,+}}(H^2(M,\R))$.
There is a $\C P^1$ of 
oriented 2-dimensional planes in $W$;
this family is precisely the twistor family
associated with the hyperk\"ahler structure
corresponding to $W$.

In the present paper, we consider what happens if one takes
a 3-di\-men\-sional plane $W\subset H^2(M,\R)$
with a degenerate metric of signature $(+,+,0)$.
Instead of a $\C P^1$ worth of complex structures,
as happens when $W$ is positive, the set of positive
2-planes in $W\subset H^2(M,\R)$ is parametrized by $\C=\R^2$.
It turns out that the corresponding family can be constructed
explicitly from an appropriate semipositive form 
on a manifold, whenever such a form exists. Moreover, this family (called 
{\bf a degenerate twistor family}; see \ref{_dege_twi_Definition_})
is holomorphic and has a canonical smooth 
trivialization, just as the usual twistor family.

\subsection{Semipositive (1,1)-forms, degenerate twistor families
and SYZ conjecture}

Let $(M,I,\Omega)$ be a simple holomorphically symplectic
manifold of K\"ahler type (that is, a hyperk\"ahler
manifold),
and $\eta\in \Lambda^{1,1}(M,I)$ a real, positive, closed
$(1,1)$-form. By Fujiki formula, either $\eta$ 
is strictly positive somewhere, or at least half of
the eigenvalues of $\eta$ vanish 
(\ref{_semipo_rank_Proposition_}). In the latter case,
the form $\Omega_t:=\Omega+t\eta$ is non-degenerate and
satisfies the assumption $\Omega_t^{n+1}=0$
for all $t$, hence defines a complex structure
(\ref{_semipis_integra_Theorem_}).

This is used to define the degenerate twistor
space (\ref{_dege_twi_integra_Theorem_}).

Positive, closed forms $\eta\in \Lambda^{1,1}(M)$
with $\int_M\eta^{\dim_\C M}=0$ are called 
{\bf semipositive}. Such forms necessarily 
lie in the boundary of a K\"ahler cone;
this implies that their cohomology classes are nef
(\ref{_nef_Definition_}). 

Notice that we exclude strictly positive forms from this definition.

\hfill

\remark
The conventions for positivity of differential forms and
currents are intrinsically confusing. Following the French tradition,
one says ``positive form'' meaning really ``non-negative'',
and ``strictly positive'' meaning ``positive definite''.
On top of it, for $(n-k, n-k)$ forms on $n$-manifold,
with $2\leq k \leq n-2$, there are two notions of 
positive forms, called ``strongly positive'' and 
``weakly positive''; this creates monsters such that
``stricly weakly positive'' and ``non-strictly stronly positive''.
The various notions of positivity in this paper
are taken (mostly) from \cite{_Demailly:Trento_},
following the French conventions as explained.

\hfill

The study of nef classes which satisfy $\int_M\eta^{\dim_\C M}=0$
(such classes are called {\bf parabolic})
is one of the central themes of hyperk\"ahler geometry.
One of the most important conjectures in this direction
is the so-called hyperk\"ahler SYZ conjecture,
due to Tyurin-Bogomolov-Hassett-Tschinkel-Huybrechts-\-Sawon
(\cite{_Hassett_Tschinkel:SYZ_conj_}, \cite{_Sawon_},
\cite{_Huybrechts:lec_}; for more history, please
see \cite{_Verbitsky:SYZ_}). This conjecture postulates
that any rational nef class $\eta$ on a hyperk\"ahler manifold
is semiample, that is, associated with a holomorphic map
$\phi:\; M \arrow X$, $\eta=\phi^*\omega_X$, where $\omega_X$
is a K\"ahler class on $X$. For nef classes which satisfy
$\int_M\eta^{\dim_\C M}>0$ (such nef classes are known
as {\bf big}), semiampleness follows from the Kawamata
base point free theorem (\cite{_Kawamata:Pluricanonical_}),
but for parabolic classes it is quite non-trivial.

If a parabolic class $\eta$  is semiample, 
it can obviously be represented by a 
smooth, semipositive differential form.
The converse implication is not proven. However,
in \cite{_Verbitsky:SYZ_} it was shown that
whenever a rational parabolic class can be
represented by a semipositive form, it is
$\Q$-effective (that is, represented by a
rational effective divisor).

Existence of a smooth semipositive form in a given
nef class is a separate (and interesting) question of
hyperk\"ahler geometry. The following conjecture is supported
by empirical evidence obtained by S. Cantat and Dinh-Sibony
(\cite{_Cantat:Acta-2001_}, 
\cite[Theorem 5.3]{_Cantat:Milnor-survey_}, 
\cite[Corollary 3.5]{_Dinh_Sibony:Jams_2005_}).

\hfill

\conjecture
Let $\eta$ be a parabolic nef class on a hyperk\"ahler manifold. 
Then $\eta$ can be represented by a semipositive closed form
with mild (say, H\"older) singularities.

\hfill

Notice that $\eta$ can be represented by a closed, positive
current by compactness of the space of positive currents
with bounded mass; however, there is no clear way to 
understand the singularities of this current.

If this conjecture is true, a 
cohomology class is $\Q$-effective whenever it is
nef and rational (\cite{_Verbitsky:SYZ_}, 
\cite{_Verbitsky:parabolic_}); 
this would prove a part of SYZ conjecture.

One of the ways of representing a nef class
by a semipositive form is based on reverse-engineering
the construction of degenerate twistor spaces. 
Let $\eta$ be a parabolic nef class on a hyperk\"ahler
manifold $(M,I)$, $\Omega$ its holomorphic symplectic
form, and $W:=\langle \eta,\Re\Omega, \Im\Omega\rangle$
the corresponding 3-dimensional subspace in $H^2(M,\R)$.
Clearly, the Bogomolov-Beauville-Fujiki form  on $W$ is degenerate
of signature $(+,+,0)$. The set $S$ of positive, oriented 
2-dimensional planes $V\subset W$ is parametrized by $\C$. Identifying
the Grassmannian $\Gr_{++}(H^2(M,\R))$ with 
a component of $\Teich_b$  as in \ref{_hk_Torelli_Theorem_},
we obtain a deformation ${\cal Z}\arrow S$; as explained in
Subsection \ref{_dege_twi_Teich_Subsection_}, 
this family can be obtained as a limit of twistor
families. The twistor families are split as smooth manifolds:
$\Tw(M)=M\times \C P^1$; this gives an Ehresmann connection $\nabla$ on 
the twistor family $\Tw(M)\arrow \C P^1$. This connection
satisfies $\nabla\Omega_t=\lambda\omega_I$, that is, 
a derivative of a holomorphically symplectic form
is proportional to a K\"ahler form. If this connection 
converges to a smooth connection $\nabla_0$ on the 
limit family ${\cal Z}\arrow \C$, we would
obtain $\nabla\Omega_t=\lambda\eta$, where $\eta$
is a limit of K\"ahler forms, hence semipositive.
This was the original motivation for the study
of degenerate twistor spaces.

\subsection{Degenerate twistor spaces and Lagrangian fibrations}

The main source of examples of degenerate twistor families comes from
Lagrangian fibrations.

\hfill

Let $(M,\Omega)$ be a simple holomorphically symplectic
K\"ahler manifold, and $\phi:\; M \arrow X$ a surjective holomorphic
map, with $0<\dim X < \dim M$. Matsushita 
(\ref{_Matsushita_fibra_Theorem_}) has shown that 
$\phi$ is a Lagrangian fibration, that is, the
fibers of $\phi$ are Lagrangian subvarieties in $M$,
and all smooth fibers of $\phi$ are Lagrangian tori. 
It is not hard to see that $X$ is projective
(\cite{_Matsushita:CP^n_}). Let $\omega_X$
be the K\"ahler form on $X$. Then 
$\eta:=\phi^*\omega_X$ is a semipositive form,
and \ref{_semipis_integra_Theorem_}
together with \ref{_comple_stru_from_2-form_Theorem_}
imply existence of a degenerate twistor family
${\cal Z}\arrow \C$, with the fibers 
holomorphically symplectic manifolds
$(M,\Omega+t\eta)$, $t\in \C$.
For each fiber $Y:=\phi^{-1}(y)$,
the restriction $\eta\restrict Y$
vanishes, because $\eta=\phi^*\omega_X$.
Therefore, the complex structure induced
by $\Omega_t=\Omega+t\eta$
on $Y$ does not depend on $t$.
This implies that the fibers of $\phi$ 
remain holomorphic and independent from $t\in \C$.

\hfill

\theorem\label{_lagra_fibra_independent_Theorem_}
Let $M$ be a simple hyperk\"ahler manifold
equipped with a Lagrangian fibration
$\phi:\; M\arrow X$, and $(M_t,\Omega_t)$
the degenerate twistor deformation associated
with the family of non-degenerate 2-forms
$\Omega+t\eta$, $\eta=\phi^*\omega_X$ as 
in \ref{_semipis_integra_Theorem_}.
Then the fibration $M_t\stackrel{\phi_t}\arrow X$
is also holomorphic, 
and for any fixed $x\in X$,  the fibers of $\phi_t$ are naturally
isomorphic: $\phi_t^{-1}(x)\cong\phi^{-1}(x)$ for all $t\in \C$.

\hfill

{\bf Proof:} The complex structure on $M_t$ is determined
from $T^{0,1}M_t=\ker \Omega_t$. Let $Z:=\phi^{-1}(x)$.
Since $\eta(v,\cdot)=0$ for each $v\in T_zZ$, 
one has $TZ\cap \ker \Omega_t=T^{0,1}Z$,
hence the complex structure on $Z$ is independent from $t$.
Since $Z$ is Lagrangian in $M_t$, its normal bundle is
dual to $TZ$ and trivial when $Z$ is a torus (that is,
for all smooth fibers of $\phi$). Therefore, the 
complex structure on $NZ$ 
is independent from $t\in \C$. This implies that the
 projection $M_t\stackrel \phi\arrow X$
is holomorphic in the smooth locus of $\phi$
for all $t\in \C$. To extend it to the points
where $\phi$ is singular, we notice that a
map is holomorphic whenever its differential
is complex linear, and complex linearity of a given tensor needs
to be checked only in an open dense subset.
\endproof

\hfill

\remark
In \cite{_Markman:Lag-2013_}, Eyal Markman
considered the following procedure.
One starts with a Lagrangian fibration $\pi$ on a hyperk\"ahler
manifold and takes a 1-cocycle on the base of $\pi$ taking values in
fiberwise automorphisms of the fibration. Twisting the $\pi$
by such a cocycle, one obtains another  Lagrangian fibration
with the same base and the respective fibers isomorphic to that
of $\pi$. Markman calls this procedure 
``the Tate-Shafarevich twist''. In this context, degenerate twistor
deformations associated with semipositive forms $\eta$,
$[\eta]\in H^2(M,\Z)$, occur very naturally; Markman calls them 
``Tate-Shafarevich lines''. One can view $\eta=\phi^*\omega_X$
as lying in
\[
\phi^*H^{1,1}(X)= \phi^* H^1(X, \Omega^1X)\subset
H^1(M, \phi^* \Omega^1X)= H^1(M, T_{M/X}),
\]
where $T_{M/X}$ is the fiberwise tangent bundle, and 
$\phi^* \Omega^1X=T_{M/X}$ because $M \arrow X$
is a Lagrangian fibration. Of course, this cocycle comes from 
$X$ so it is constant in the fibre direction; it
describes the deformation infinitesimally. Integrating 
the vector field then gives a 1-cocycle on
$X$ taking values in the bundle of fibrewise 
automorphisms. This is the 1-cocycle giving the
"Tate-Shafarevich twist".

\hfill

\remark
The degenerate twistor family constructed in \ref{_dege_twi_integra_Theorem_} 
consists of a family of complex structures, but it is not proven that all fibers,
which are complex manifolds, are also K\"ahler (hence hyper\"ahler).
As is, the K\"ahler property is known only over a small open subset
in the base (affine line), since the condition of being K\"ahler is open.
We expect all members of the degenerate twistor family to be
K\"ahler, but there is no obvious way to prove this. However,
it is easy to show that the set of points on the base affine line
corresponding to non-K\"ahler complex structures is 
closed and countable.


\section{Basic notions of hyperk\"ahler geometry}


\subsection{Hyperk\"ahler manifolds}
\label{_hk_basic_Subsection_}

\definition
Let $(M,g)$ be a Riemannian manifold, and $I,J,K$
endomorphisms of the tangent bundle $TM$ satisfying the
quaternionic relations
\[
I^2=J^2=K^2=IJK=-\Id_{TM}.
\]
The triple $(I,J,K)$ together with
the metric $g$ is called {\bf a hyperk\"ahler structure}
if $I, J$ and $K$ are integrable and K\"ahler with respect to $g$.

Consider the K\"ahler forms $\omega_I, \omega_J, \omega_K$
on $M$:
\begin{equation}\label{_omega_I,J,K_defi_Equation_}
\omega_I(\cdot, \cdot):= g(\cdot, I\cdot), \ \
\omega_J(\cdot, \cdot):= g(\cdot, J\cdot), \ \
\omega_K(\cdot, \cdot):= g(\cdot, K\cdot).
\end{equation}
An elementary linear-algebraic calculation implies
that the 2-form 
\begin{equation}\label{_holo_symple_on_hk_Equation_}
\Omega:=\omega_J+\1\omega_K
\end{equation}
is of Hodge type $(2,0)$
on $(M,I)$. This form is clearly closed and
non-degenerate, hence it is a holomorphic symplectic form.

In algebraic geometry, the word ``hyperk\"ahler''
is essentially synonymous with ``holomorphically
symplectic'', due to the following theorem, which is
implied by Yau's solution of Calabi conjecture
(\cite{_Beauville_,_Besse:Einst_Manifo_}).

\hfill

\theorem\label{_Calabi-Yau_Theorem_}
Let $M$ be a compact, K\"ahler, holomorphically
symplectic manifold, $\omega$ its K\"ahler form, $\dim_\C M =2n$.
Denote by $\Omega$ the holomorphic symplectic form on $M$.
Assume that $\int_M \omega^{2n}=\int_M (\Re\Omega)^{2n}$.
Then there exists a unique hyperk\"ahler metric $g$ within the same
K\"ahler class as $\omega$, and a unique hyperk\"ahler structure
$(I,J,K,g)$, with $\omega_J = \Re\Omega$, $\omega_K = \im\Omega$.
\endproof 

\hfill

\subsection{The Bogomolov-Beauville-Fujiki form}

\definition\label{_hk_simple_Definition_}
A hyperk\"ahler manifold $M$ is called
{\bf simple} if $\pi_1(M)=0$, $H^{2,0}(M)=\C$.
In the literature, such manifolds are often
called {\bf irreducible holomorphic symplectic},
or {\bf irreducible symplectic varieties}.

\hfill

This definition is motivated by the following theorem
of Bogomolov (\cite{_Bogomolov:decompo_}). 

\hfill

\theorem (\cite{_Bogomolov:decompo_})
Any hyperk\"ahler manifold admits a finite covering
which is a product of a torus and several 
simple hyperk\"ahler manifolds.
\endproof

\hfill

\theorem\label{_Fujiki_Theorem_}
(\cite{_Fujiki:HK_})
Let $\eta\in H^2(M)$, and $\dim M=2n$, where $M$ is
a simple hyperk\"ahler manifold. Then $\int_M \eta^{2n}=\lambda q(\eta,\eta)^n$,
for some integer quadratic form $q$ on $H^2(M)$, and $\lambda\in\Q$
a positive rational number.
\endproof

\hfill

\definition\label{_BBF_Definition_}
This form is called
{\bf  Bogomolov-Beauville-Fujiki form}.  It is defined
by this relation uniquely, up to a sign. The sign is determined
from the following formula (Bogomolov, Beauville;
\cite{_Beauville_}, \cite{_Huybrechts:lec_}, 23.5)
\begin{align*}   \lambda q(\eta,\eta) &=
   (n/2)\int_X \eta\wedge\eta  \wedge \Omega^{n-1}
   \wedge \bar \Omega^{n-1} -\\
 &-(1-n)\frac{\left(\int_X \eta \wedge \Omega^{n-1}\wedge \bar
   \Omega^{n}\right) \left(\int_X \eta \wedge
   \Omega^{n}\wedge \bar \Omega^{n-1}\right)}{\int_M \Omega^{n}
   \wedge \bar \Omega^{n}}
\end{align*}
where $\Omega$ is the holomorphic symplectic form,
and $\lambda$ a positive constant.

\hfill

\remark
The form $q$ has signature $(3,b_2-3)$.
It is negative definite on primitive forms, and positive
definite on the space $\langle \Re \Omega, \Im\Omega, \omega\rangle$
 where $\omega$ is a K\"ahler form, as seen from the
following formula
\begin{multline}\label{_BBF_via_Kahler_Equation_}
   \mu q(\eta_1,\eta_2)= \\
   \int_X \omega^{2n-2}\wedge \eta_1\wedge\eta_2  
   - \frac{2n-2}{(2n-1)^2}
   \frac{\int_X \omega^{2n-1}\wedge\eta_1 \cdot
   \int_X\omega^{2n-1}\wedge\eta_2}{\int_M\omega^{2n}}, \
   \  \mu>0
\end{multline}
(see e. g. \cite{_Verbitsky:cohomo_}, Theorem 6.1,
or \cite{_Huybrechts:lec_}, Corollary 23.9).

\hfill

\definition
Let $[\eta]\in H^{1,1}(M)$ be a real (1,1)-class in
the closure of the K\"ahler cone of
a hyperk\"ahler manifold $M$. We say that $[\eta]$
is {\bf parabolic} if $q([\eta],[\eta])=0$.

\subsection{The hyperk\"ahler SYZ conjecture}

\theorem\label{_Matsushita_fibra_Theorem_}
(D. Matsushita, see \cite{_Matsushita:fibred_}).
Let $\pi:\; M \arrow X$ be a surjective holomorphic map
from a simple hyperk\"ahler manifold $M$ to 
a complex variety $X$, with $0<\dim X < \dim M$.
Then $\dim X = 1/2 \dim M$, and the fibers of $\pi$ are 
holomorphic Lagrangian (this means that the symplectic
form vanishes on the fibers).\footnote{Here, as elsewhere,
we silently assume that the hyperk\"ahler manifold $M$ is
simple.}

\hfill

\definition Such a map is called
{\bf a holomorphic Lagrangian fibration}.

\hfill

\remark The base of $\pi$ is conjectured to be
rational. J.-M. Hwang (\cite{_Hwang:base_}) 
proved that $X\cong \C P^n$, if $X$ is smooth
and $M$ projective.
D. Matsushita (\cite{_Matsushita:CP^n_}) 
proved that it has the same rational cohomology
as $\C P^n$ when $M$ is projective.

\hfill

\remark
 The base of $\pi$ has a natural flat connection
on the smooth locus of $\pi$. The combinatorics of this connection
can be (conjecturally) used to determine the topology of $M$ 
(\cite{_Kontsevich-Soibelman:torus_},  
\cite{_Kontsevich-Soibelman:non-archimedean_},
\cite{_Gross:SYZ_}).

\hfill

\remark 
Matsushita's theorem is implied by the following formula
of Fujiki. Let $M$ be a hyperk\"ahler manifold, $\dim_\C
M=2n$, and $\eta_1, ..., \eta_{2n}\in H^2(M)$ cohomology
classes. Then 
\begin{equation}\label{_Fujiki_multi_Equation_}
C\int_M \eta_1\wedge \eta_2 \wedge ... = 
\frac{1}{(2n)!}\sum_{\sigma}q(\eta_{\sigma_1} \eta_{\sigma_2})
q(\eta_{\sigma_3}\eta_{\sigma_4})... q(\eta_{\sigma_{2n-1}} \eta_{\sigma_{2n}})
\end{equation}
with the sum taken over all permutations, and $C$
a positive constant, called {\bf Fujiki constant}.
An algebraic argument (see e.g. \ref{_product_vanishes_Corollary_})
allows to deduce from this formula that
for any non-zero $\eta \in H^2(M)$,
one would have  $\eta^{n}\neq 0$, and $\eta^{n+1}=0$, 
if $q(\eta,\eta)=0$, and $\eta^{2n}\neq 0$ otherwise.
Applying this to the pullback $\pi^*\omega_X$
of the K\"ahler class from $X$, we immediately
obtain that $\dim_\C X=n$ or $\dim_\C X=2n$. Indeed, 
$\omega_X^{\dim_\C X}\neq 0$ and
$\omega_X^{\dim_\C X+1}= 0$. This argument
was used by Matsushita in his proof of 
\ref{_Matsushita_fibra_Theorem_}.
The relation \eqref{_Fujiki_multi_Equation_}
is another form of Fujiki's theorem (\ref{_Fujiki_Theorem_}),
obtained by differentiation of 
$\int_M \eta^{2n}=\lambda q(\eta,\eta)^n$,

\subsection{Cohomology of hyperk\"ahler manifolds}

Further on in this paper, some basic results about cohomology
of hyperk\"ahler manifolds will be used.
The following theorem was proved in \cite{_Verbitsky:cohomo_},
using representation theory.

\hfill

\theorem \label{_symme_coho_Theorem_}
(\cite{_Verbitsky:cohomo_}) 
Let $M$ be a simple hyperk\"ahler manifold, 
and $H^*_r(M)$ the part of cohomology generated by $H^2(M)$.
Then $H^*_r(M)$ is isomorphic to the symmetric algebra 
(up to the middle degree). Moreover, the Poincare 
pairing on $H^*_r(M)$ is non-degenerate. 
\endproof

\hfill
 
This brings the following corollary.

\hfill

\corollary \label{_product_vanishes_Corollary_}
Let $\eta_1, ... \eta_{n+1}\in H^2(M)$ be cohomology
classes on a simple hyperk\"ahler manifold, $\dim_\C M = 2n$.
Suppose that $q(\eta_i, \eta_j)=0$ for all $i, j$. Then
$\eta_1 \wedge \eta_2 \wedge ... \wedge \eta_{n+1}=0$.

{\bf Proof:} See e.g. \cite[Corollary 2.15]{_Verbitsky:parabolic_}. 
This equation also follows from \eqref{_Fujiki_multi_Equation_}.
\endproof

\section{Degenerate twistor space}

\subsection{Integrability of almost complex structures and Cartan formula}
\label{_integra_Subsection_}

An {\bf almost complex structure} on a manifold is
a section $I\in \End(TM)$ of the bundle of endomorphisms,
satisfying $I^2=-\Id$. It is called {\bf integrable}
if $[T^{1,0}M,T^{1,0}M]\subset T^{1,0}M$, where
$T^{1,0}M\subset TM\otimes_\R\C$ is the eigenspace of $I$,
defined by 
\[
v\in T^{1,0}M \Leftrightarrow I(v)=\1 v.
\]
Equivalently, $I$ is integrable if 
$[T^{0,1}M,T^{0,1}M]\subset T^{0,1}M$,
where $T^{0,1}M\subset TM\otimes_\R\C$ 
is a complex conjugate to $T^{1,0}M\subset TM\otimes_\R\C$.

One of the ways of making sure a given almost complex structure
is integrable is by using the Cartan formula expressing the de Rham
differential through commutators of vector fields.

\hfill

\proposition\label{_2-form_closed_the_integra_Proposition_}
Let $(M,I)$ be a manifold equipped 
with an almost complex structure, and $\Omega\in \Lambda^{2,0}(M)$ 
a non-degenerate $(2,0)$-form (\ref{_non-dege_Definition_}). 
Assume that $d\Omega=0$.
Then $I$ is integrable.

\hfill

{\bf Proof:} Let $X \in T^{1,0}M$ and $Y,Z\in T^{0,1}(M)$.
Since $\Omega$ is a (2,0)-form, it vanishes on 
$(0,1)$-vectors. Then Cartan formula together with 
$d\Omega=0$ implies that
\begin{equation}\label{_Cartan_Eqution_}
0=d\Omega(X,Y,Z)= \Omega(X,[Y,Z]).
\end{equation}
From the non-degeneracy of $\Omega$ we obtain that
unless $[Y,Z]\in T^{0,1}(M)$,
for some $X \in T^{1,0}M$,
one would have $\Omega(X,[Y,Z])\neq 0$.
Therefore, \eqref{_Cartan_Eqution_} implies
 $[Y,Z]\in T^{0,1}(M)$, for all $Y,Z\in T^{0,1}(M)$,
which means that $I$ is integrable. 
\endproof

\hfill

\remark\label{_2,1_diffe_symple_Remark_}
It is remarkable that the closedness of
$\Omega$ is in fact unnecessary.
The proof \ref{_2-form_closed_the_integra_Proposition_}
remains true if one assumes that
$d\Omega\in \Lambda^{3,0}(M)\oplus \Lambda^{2,1}(M)$.

\hfill

Notice that the sub-bundle $T^{1,0}M\subset
TM\otimes_\R\C$ uniquely determines the almost complex
structure. Indeed, $I(x+ y)=\1 x-\1 y$, for all
$x\in T^{1,0}M, y \in T^{0,1}M=\overline{T^{1,0}M}$,
and we have a decomposition $T^{1,0}M\oplus T^{0,1}M=TM\otimes_\R\C$.
This decomposition is the necessarily and sufficient
ingredient for the reconstruction of an almost complex
structure:

\hfill

\claim\label{_almost_co_by_(0,1)_Claim_}
Let $M$ be a smooth, $2n$-dimensional manifold. Then
there is a bijective correspondence between the set
of almost complex structures, and the set of sub-bundles
$T^{0,1}M\subset TM\otimes_\R\C$ satisfying $\dim_\C T^{0,1}M= n$
and $T^{0,1}M\cap TM=0$ (the last condition means that
there are no real vectors in $T^{1,0}M$).
\endproof

\hfill

The last two statements allow us to define complex structures
in terms of complex-valued 2-forms (see
\ref{_comple_stru_from_2-form_Theorem_} below).
For this theorem, any reasonable notion of non-degeneracy would
suffice; for the sake of clarity, we state the one we would use.

\hfill

\definition\label{_non-dege_Definition_}
Let  $\Omega\in \Lambda^2(M,\C)$ be a smooth, complex-valued 2-form on a 
$2n$-dimensional manifold. $\Omega$ is called {\bf non-degenerate}
if for any real vector $v\in T_mM$, the contraction $\Omega\cntrct v$
is non-zero.

\hfill

\theorem\label{_comple_stru_from_2-form_Theorem_}
Let $\Omega\in \Lambda^2(M,\C)$ be a smooth, complex-valued,
non-degenerate 2-form on a $4n$-dimensional real manifold. 
Assume that $\Omega^{n+1}=0$. Consider the bundle
\[
T^{0,1}_\Omega(M):= \{ v\in TM\otimes \C \ \ |\ \  \Omega\cntrct v=0\}.
\]
Then $T^{0,1}_\Omega(M)$ satisfies assumptions of
\ref{_almost_co_by_(0,1)_Claim_}, hence defines an 
almost complex structure $I_\Omega$ on $M$. If, in addition,
$\Omega$ is closed, $I_\Omega$ is integrable.

\hfill

{\bf Proof:}
Integrability of $I_\Omega$ follows immediately
from \ref{_2-form_closed_the_integra_Proposition_}.
Let $v\in TM$ be a non-zero real tangent vector. Then
$\Omega\cntrct v\neq 0$, hence $T^{0,1}_\Omega(M)\cap TM=0$.
To prove \ref{_comple_stru_from_2-form_Theorem_},
it remains to show that $\rk T^{0,1}_\Omega(M)\geq 2n$.
Clearly, $\Omega$ is non-degenerate on $\frac{TM\otimes \C}{T^{0,1}_\Omega(M)}$,
hence its rank is equal to $4n-\rk T^{0,1}_\Omega(M)$.
From $\Omega^{n+1}=0$ it follows that rank of $\Omega$ cannot
exceed $2n$, hence $\rk T^{0,1}_\Omega(M)\geq 2n$.
\endproof

\subsection{Semipositive (1,1)-forms on hyperk\"ahler manifold}

\definition
Let $\eta\in \Lambda^{1,1}(M,\R)$ be a real (1,1)-form on a complex
manifold $(M,I)$.
It is called {\bf semipositive} if $\eta(x,Ix)\geq 0$ for any
$x\in TM$, but it is nowhere positive definite. 

\hfill

\remark
Fix a Hermitian structure $h$ on $(M,I)$.
Clearly, any semipositive (1,1)-form is diagonal in some
$h$-orthonormal basis in $TM$. The entries of its matrix 
in this basis are called {\bf eigenvalues}; they are
real, non-negative numbers. The maximal number of 
positive eigenvalues is called {\bf the rank} of a 
semipositive (1,1)-form.

\hfill

\definition\label{_nef_Definition_}
A closed semipositive form $\eta$ on a compact K\"ahler manifold $(M,I,\omega)$
is a limit of K\"ahler forms $\eta+\epsilon \omega$, hence
its cohomology class is {\bf nef} (belongs to the closure of the
K\"ahler cone). Its cohomology class $[\eta]$ is {\bf parabolic},
that is, satisfies $\int_M[\eta]^{\dim_\C M}=0$.
However, not every parabolic nef class can be represented
by a closed semipositive form (\cite{_Demailly_Peternell_Schneider:nef_}).

\hfill

\proposition\label{_semipo_rank_Proposition_}
On a simple
 hyperk\"ahler manifold $M$, $\dim_\C M=2n$, any semipositive (1,1)-form
has rank $0$ or $2n$.

\hfill

{\bf Proof:} This assertion
easily follows from \ref{_product_vanishes_Corollary_}.
Indeed, if $q(\eta,\eta)\neq 0$, one has 
$\int_M \eta^{2n}=\lambda q(\eta,\eta)^n\neq 0$,
hence its rank is $4n$. If $q(\eta,\eta)=0$,
its cohomology class $[\eta]$ satisfies $[\eta]^n\neq 0$ and
$[\eta]^{n+1}=0$ (\ref{_product_vanishes_Corollary_}). Since
all eigenvalues of $\eta$ are non-negative, 
its rank is twice the biggest number $k$ for which one
has $\eta^k\neq 0$. However, since $\eta^k$ is a sum
of monomials of an orthonormal basis with non-negative 
coefficients, $\int_M \eta^k\wedge \omega^{2n-k}=0$ 
$\Leftrightarrow$ $\eta^k=0$ for any K\"ahler form
$\omega$ on $(M,I)$. Then $[\eta]^n\neq 0$ and
$[\eta]^{n+1}=0$ imply that the rank of $\eta$ is $2n$.
\endproof

\hfill

The main technical result of this paper is the following theorem.

\hfill

\theorem\label{_semipis_integra_Theorem_}
Let $(M,\Omega)$ be an simple hyperk\"ahler manifold, $\dim_\R M=4n$, 
and $\eta\in \Lambda^{1,1}(M,I)$ a closed, semipositive form
of rank $2n$. Then the 2-form $\Omega+t\eta$ satisfies
the assumptions of \ref{_comple_stru_from_2-form_Theorem_}
for all $t\in \C$: namely, $\Omega+t\eta$ is non-degenerate,
and $(\Omega+t\eta)^{n+1}=0$.

\hfill

{\bf Proof:}
Non-degeneracy of $\Omega_t:=\Omega+t\eta$ is clear.
Indeed, let $v:= |t|t^{-1}$, and let 
$\omega_v:=\Re v \omega_K-\im v \omega_J$. Then  
$\omega_v$ is a Hermitian form associated with the 
induced complex structure $\Im v J- \Re v K$, hence
it is non-degenerate. However,
the imaginary part of $v\Omega_t$ is equal to 
$\omega_v$ (see \eqref{_omega_I,J,K_defi_Equation_}).
Then $\Im(\Omega_t\cntrct v)\neq 0$ for each non-zero
real vector $v\in TM$.

To see that $(\Omega+t\eta)^{n+1}=0$, we observe that
this relation is true in cohomology; this is
implied from \cite{_Verbitsky:coho_announce_} using
the same argument as was used in the proof 
of \ref{_semipo_rank_Proposition_}.

Each Hodge component of $(\Omega+t\eta)^{n+1}$ is proportional
to $\Omega^{n-p}\wedge \eta^{p+1}$, and it is sufficient
to prove that $\Omega^{n-p}\wedge \eta^{p+1}=0$ for all $p$.

We deduce this from two observations, which are proved 
further on in this section.

\hfill

\lemma\label{_posi_then_ome_vani_Lemma_}
Let $(M,\Omega)$, $\dim_\R M=4n$ be a holomorphically symplectic manifold, 
and $\eta\in \Lambda^{1,1}(M,I)$ a closed, semipositive form
of rank $2n$. Assume that $\Omega^{n-p}\wedge \eta^{p+1}$ is exact.
Then \[ \Omega^{n-p}\wedge \bar \Omega^{n-p}\wedge\eta^{p+1}=0,\]
for all $p$.

{\bf Proof:} See Subsection \ref{_posi_forms_Subsection_}. \endproof

\hfill

\lemma\label{_ome_bar_ome_vani__then_ome_vani_Lemma_}
Let $(M,\Omega)$, $\dim_\R M=4n$, be a holomorphically symplectic manifold
and $\rho\in \Lambda^{p+1,p+1}(M,I)$ a strongly positive
form (\ref{_posi_form_Definition_}). 
Suppose that $\Omega^{n-p}\wedge \bar \Omega^{n-p}\wedge \rho=0$.
Then $\Omega^{n-p}\wedge \rho=0$.

{\bf Proof:} See Subsection \ref{_posi_forms_wedge_Omega_Subsection_}. 
\endproof

\subsection{Positive $(p,p)$-forms}
\label{_posi_forms_Subsection_}

We recall the definition of a positive $(p,p)$-form
(see e.g. \cite{_Demailly:Trento_}).

\hfill

\definition\label{_posi_form_Definition_}
Recall that a real $(p,p)$-form $\eta$
on a complex manifold is called {\bf weakly positive}
if for any complex subspace $V\subset T M$, 
$\dim_\C V=p$, the restriction $\rho\restrict V$
is a non-negative volume form. Equivalently,
this means that 
\[ 
  (\1)^p\rho(x_1, \bar x_1, x_2, \bar x_2, ..., x_p, \bar
  x_p)\geq 0,
\]
for any vectors $x_1, ... x_p\in T_x^{1,0}M$.
A real $(p,p)$-form on a complex manifold 
is called {\bf strongly positive} if it can 
be locally expressed as a sum
\[
\eta = (\1)^p\sum_{i_1, ... i_p} 
\alpha_{i_1, ... i_p} \xi_{i_1} \wedge \bar\xi_{i_1}\wedge ... 
\wedge \xi_{i_p} \wedge \bar\xi_{i_p}, \ \  
\]
running over some set of $p$-tuples 
$\xi_{i_1}, \xi_{i_2}, ..., \xi_{i_p}\in \Lambda^{1,0}(M)$,
with $\alpha_{i_1, ..., i_p}$ real and non-negative
functions on $M$. 

The following basic linear algebra observations are easy to 
check (see \cite{_Demailly:Trento_}).

All strongly positive forms are also weakly positive.
The strongly positive and the weakly positive forms
form closed, convex cones in the space 
$\Lambda^{p,p}(M,\R)$ of real $(p,p)$-forms.
These two cones are dual with respect to the Poincare pairing
\[
\Lambda^{p,p}(M,\R) \times \Lambda^{n-p,n-p}(M,\R)\arrow \Lambda^{n,n}(M,\R)
\]
For (1,1)-forms and $(n-1,n-1)$-forms,
the strong positivity is equivalent
to weak positivity. Finally, a product of a weakly
positive form and a strongly positive one is always
weakly positive (however, a product of two weakly positive
forms may be not weakly positive).

\hfill

Clearly, an exact weakly positive form $\eta$ on a compact
K\"ahler manifold $(M,\omega)$ always vanishes. Indeed, 
the integral $\int_M\eta\wedge \omega^{\dim M-p}$ 
is strictly positive for a non-zero weakly positive
$\eta$, because the convex cones of weakly and strongly
positive forms are dual, and $\omega^{\dim M-p}$ sits in 
the interior of the cone of strongly positive forms.
However, by Stokes' formula, this integral vanishes 
whenever $\eta$ is exact.

Now we are in position to prove \ref{_posi_then_ome_vani_Lemma_}.
The form $\Omega^{n-p}\wedge \bar \Omega^{n-p}\wedge \eta^{p+1}$ is
by assumption of this lemma exact, but it is a product
of a weakly positive form $\Omega^{n-p}\wedge \bar \Omega^{n-p}$
and a strongly positive form $\eta^{p+1}$, hence it is weakly
positive. Being exact, this form must vanish. 

\hfill

\remark
A form is strongly positive if it is generated by products
of $dz_i\wedge d\bar z_i$ with positive coefficients;
hence $\eta$ and all its powers are positive. The form
$\Omega\wedge \bar \Omega$ and its powers are positive
on all complex spaces of appropriate dimensions, which
can be seen by using Darboux coordinates. This means
that this form is weakly positive. 

\hfill

\subsection{Positive $(p,p)$-forms and holomorphic symplectic forms}
\label{_posi_forms_wedge_Omega_Subsection_}

Now we shall prove
\ref{_ome_bar_ome_vani__then_ome_vani_Lemma_}.
This is a linear-algebraic statement, which can be proven
pointwise. Fix a complex vector space
$V$, equipped with a non-degenerate complex linear
2-form $\Omega$. Every strongly positive form $\rho$ on $V$
is a sum of monomials
$(\1)^p\xi_{i_1} \wedge \bar\xi_{i_1}\wedge ... 
\wedge \xi_{i_p} \wedge \bar\xi_{i_p}$ with positive coefficients,
and the equivalence
\[
  \Omega^{n-p}\wedge \rho\neq 0 \Leftrightarrow
\Omega^{n-p}\wedge \bar \Omega^{n-p}\wedge \rho\neq 0
\]
is implied by the following sublemma.

\hfill

\sublemma\label{_posi_mult_vector_space_Sublemma_}
Let $V$ be a complex vector space, equipped 
with a non-degenerate complex linear
2-form $\Omega\in \Lambda^{2,0}V$. Then
for any monomial
$\rho=(\1)^p\xi_{i_1} \wedge \bar\xi_{i_1}\wedge ... 
\wedge \xi_{i_p} \wedge \bar\xi_{i_p}$ for which
$\Omega^{n-p}\wedge \rho$ is non-zero, the form
$\Omega^{n-p}\wedge \bar \Omega^{n-p}\wedge \rho$
is non-zero and weakly positive. 

\hfill

{\bf Proof:}
Let $\xi_{j_1}, \xi_{j_1}, ...,
\xi_{j_{n-p}}$ be the elements of the basis in $V$
complementary to $\xi_{i_1}, \xi_{i_1}, ...,
\xi_{i_p}$, and $W\subset V$ the space generated
by $\xi_{j_1}, \xi_{j_1}, ...,
\xi_{j_{n-p}}$. Clearly, a form $\alpha$ is non-zero
on $W$ if and only if $\alpha\wedge \rho$ is non-zero,
and positive on $W$ if and only if $\alpha\wedge \rho$ is positive.

Now, \ref{_posi_mult_vector_space_Sublemma_}
is implied by the following trivial assertion:
for any $(n-p)$-dimensional 
subspace $W\subset V$ such that $\Omega^{n-p}\restrict W$
is non-zero, the restriction $\Omega^{n-p}\wedge \bar \Omega^{n-p}\restrict W$
is non-zero and positive.

This proves \ref{_posi_mult_vector_space_Sublemma_},
and \ref{_ome_bar_ome_vani__then_ome_vani_Lemma_} follows
as indicated. \endproof

\hfill

As a corollary of the vanishing of the forms $\Omega^{n-p}\wedge\eta^{p+1}$,
we prove the following statement, used further on.

\hfill

\lemma\label{_eta_1,1_Lemma}
Let $(M,\Omega)$ be a simple
holomorphically symplectic manifold, $\dim_\R M=4n$ 
and $\eta\in \Lambda^{1,1}(M,I)$ a closed, semipositive form
of rank $2n$. 
Let $I_t$ be the complex structure on $M$ defined by
$\Omega+t\eta$, as in \ref{_semipis_integra_Theorem_}. 
Then $\eta\in \Lambda^{1,1}(M,I_t)$.

\hfill

{\bf Proof:} By construction, $(M,I_t)$ is a holomorphically
symplectic manifold, with the holomorphic symplectic form
$\Omega_t:=\Omega+t\eta$. For a holomorphic symplectic manifold
$(M,\Omega_t)$, $\dim_\R M=4n$,
there exist an elementary criterion allowing one to check whether
a given 2-form $\eta$ is of type (1,1): one has to have
$\eta\wedge \Omega_t^n=0$ and $\eta\wedge \bar \Omega_t^n=0$.
However, from \ref{_ome_bar_ome_vani__then_ome_vani_Lemma_} it follows
immediately that $\eta\wedge \Omega_t^n=0$ and 
$\eta\wedge \bar\Omega_t^n=0$, hence $\eta$ is of type (1,1).
\endproof

\subsection{Degenerate twistor space: a definition}
\label{_dege_twi_Subsection_}

Just as it is done with the usual twistor space,
to define a degenerate twistor space we construct
a certain almost complex structure, and then prove
it is integrable. The proof of integrability is in fact 
identical to the argument which could be used to prove that the
usual twistor space is integrable.

\hfill

\definition\label{_dege_twi_Definition_}
Let $(M,\Omega)$ be an irreducible 
holomorphically symplectic manifold, $\dim_\R M=4n$ 
and $\eta\in \Lambda^{1,1}(M,I)$ a closed, semipositive form
of rank $2n$. Consider the product $\Tw_\eta(M):=\C \times M$,
equipped with the almost complex structure ${\cal I}$
acting on $T_t\C \oplus T_m M$ as $I_\C \oplus I_t$,
where $I_\C$ is the standard complex structure on $\C$
and $I_t$ is the complex structure recovered from the
form $\Omega+t\eta$ using \ref{_semipis_integra_Theorem_}
and \ref{_comple_stru_from_2-form_Theorem_}.
The almost complex manifold $(\Tw_\eta(M), {\cal I})$
is called {\bf a degenerate twistor space} of $M$.

\hfill

\theorem\label{_dege_twi_integra_Theorem_}
The almost complex structure on a degenerate twistor space is 
always integrable.

\hfill

{\bf Proof:}
We introduce a dummy variable $w$, and consider a product
$\Tw_\eta(M)\times \C$, equipped with the (2,0)-form
$\tilde \Omega:= \Omega+t \eta+ dt\wedge dw$.
Here, $\Omega$ is a holomorphic symplectic form on $M$
lifted to $M\times \C \times \C$, and $t$ and $w$ are 
complex coordinates
on $\C \times \C$. Clearly, $\tilde \Omega$ is a non-degenerate
(2,0)-form. From \ref{_eta_1,1_Lemma} we obtain
that $d\tilde\Omega=\eta\wedge dt\in\Lambda^{2,1}(\Tw_\eta(M)\times \C)$.
Now, \ref{_2,1_diffe_symple_Remark_} implies that
$\tilde \Omega$ defines an integrable almost complex
structure on $\Tw_\eta(M)\times \C$. However, on
$\Tw_\eta(M)\times\{w\}$ this almost complex structure
coincides with the one given by the degenerate twistor construction.
\endproof

\hfill

{\bf Acknowledgements:}
I am grateful to Eyal Markman 
and Jun-Muk Hwang for their interest and encouragement.
Thanks to Ljudmila Kamenova for her suggestions and
to the organizers of the Quiver Varieties Program 
at the Simons Center for Geometry and Physics,
Stony Brook University, where some of the research 
for this paper was performed. Also much gratitude
to the anonymous referee for important suggestions.

{\scriptsize

\noindent {\sc Misha Verbitsky\\
Laboratory of Algebraic Geometry, \\
Faculty of Mathematics, NRU HSE,\\
7 Vavilova Str. Moscow, Russia
\tt  verbit@mccme.ru},\\ also: \\
{\sc Kavli IPMU (WPI), the University of Tokyo.}
}


\begin{thebibliography}{GMP}

\bibitem[Bea]{_Beauville_} 
 Beauville, A. {\em 
Varietes K\"ahleriennes dont la premi\`ere classe de Chern est
nulle.}  J. Diff. Geom. {\bf 18}, pp. 755-782 (1983).


\bibitem[Bes]{_Besse:Einst_Manifo_} 
Besse, 
A., {\em Einstein Manifolds}, Springer-Verlag, New York (1987)


\bibitem[Bo1]{_Bogomolov:decompo_}  
Bogomolov, F. A., {\em On the decomposition of 
K\"ahler manifolds with trivial canonical class}, Math. USSR-Sb.
{\bf 22} (1974), 580-583.



\bibitem[C1]{_Cantat:Acta-2001_}
Serge Cantat, {\em Dynamique des automorphismes des surfaces K3}, 
Acta Mathematica,
    187:1-57, 2001.

\bibitem[C2]{_Cantat:Milnor-survey_}
Serge Cantat, {\em Dynamics of automorphisms of compact complex surfaces}, 
in "Frontiers in Complex Dynamics: a volume in honor of John Milnor's 
80th birthday", Princeton University Press.

\bibitem[D]{_Demailly:Trento_}
Demailly, Jean-Pierre,
{\em Monge-Amp\`ere operators, Lelong numbers and intersection
theory,} Complex Analysis and Geometry, Univ. Series in
Math., edited by V. Ancona and A. Silva, Plenum Press,
New-York (1993) 



\bibitem[DPS]{_Demailly_Peternell_Schneider:nef_}
Jean-Pierre Demailly, 
Thomas Peternell, Michael Schneider,
{\em Compact complex manifolds with numerically effective
tangent bundles,} J. Algebraic Geometry 3 (1994) 295-345

\bibitem[DS]{_Dinh_Sibony:Jams_2005_}
 Tien-Cuong Dinh, Nessim Sibony, 
{\em Green currents for holomorphic 
automorphisms of compact K\"ahler manifolds},
 J. Amer. Math. Soc. 18 (2005), no. 2, 291-312

\bibitem[F]{_Fujiki:HK_}  
Fujiki, A. {\em On the de Rham Cohomology Group of a Compact 
K\"ahler Symplectic Manifold}, Adv. Stud.
Pure Math. 10 (1987), 105-165.

\bibitem[G]{_Gross:SYZ_} 
Mark Gross,
{\em Mirror Symmetry and the Strominger-Yau-Zaslow conjecture}, 
arXiv:1212.4220, 67 pages.


\bibitem[HT]{_Hassett_Tschinkel:SYZ_conj_}
Brendan Hassett, Yuri Tschinkel, {\em 
Rational curves on holomorphic symplectic fourfolds},
arXiv:math/9910021,  
Geom. Funct. Anal. 11 (2001), no. 6, 1201--1228.


\bibitem[H1]{_Huybrechts:basic_} 
Huybrechts, D., 
{\em Compact hyperk\"ahler manifolds: Basic
results}, Invent. Math. 135 (1999), 63-113, 
alg-geom/9705025


\bibitem[Hu2]{_Huybrechts:lec_}
 Huybrechts, Daniel, 
{\em Compact hyperk\"ahler manifolds, Calabi-Yau 
manifolds and related geometries,}
Universitext, Springer-Verlag, Berlin, 2003, 
Lectures from the Summer School held in
Nordfjordeid, June 2001, pp. 161-225.


\bibitem[Hw]{_Hwang:base_}
Jun-Muk Hwang,
{\em Base manifolds for fibrations of projective irreducible
symplectic manifolds}, arXiv:0711.3224,
Inventiones mathematicae, vol. 174, issue 3, pp. 625-644.



\bibitem[Kal]{_Kaledin:twistor_}
D. Kaledin,
{\em Integrability of the twistor space for a hypercomplex manifold},
Selecta Math. (N.S.) {\bf 4} (1998) 271-278.


\bibitem[Kaw]{_Kawamata:Pluricanonical_}
Y. Kawamata, 
{\em Pluricanonical sysmtems on minimal algebraic varieties,} 
Invent. Math., 79, (1985), no. 3, 567-588.



\bibitem[KZ1]{_Kontsevich-Soibelman:torus_}
Maxim Kontsevich, Yan Soibelman,
{\em Homological mirror symmetry and torus fibrations},
arXiv:math/0011041, Symplectic geometry and mirror
symmetry (Seoul, 2000), World Sci. Publishing, 
River Edge, NJ, 2001, pp. 203-263.



\bibitem[KZ2]{_Kontsevich-Soibelman:non-archimedean_}
Maxim Kontsevich, Yan Soibelman,
{\em Affine structures and non-Archimedean analytic spaces}, 
The unity of mathematics, 321-385,
Progr. Math., 244, Birkh\"auser Boston, Boston, MA, 2006;
arXiv:math/0406564.

\bibitem[Mar]{_Markman:Lag-2013_}
Eyal Markman,
{\em 
Lagrangian fibrations of holomorphic-symplectic varieties of K3${}^{[n]}$-type},
arXiv:1301.6584, 34 pages.

\bibitem[Mat1]{_Matsushita:fibred_} 
D. Matsushita, {\em On fibre space structures of a
  projective irreducible symplectic manifold},
alg-geom/9709033, math.AG/9903045, also in Topology
{\bf 38} (1999), No. 1, 79-83. Addendum, 
Topology {\bf 40} (2001), No. 2, 431-432.


\bibitem[Mat2]{_Matsushita:CP^n_} 
 Matsushita, D., {\em
Higher direct images of Lagrangian fibrations},
Amer. J. Math. 127 (2005), arXiv:math/0010283.




\bibitem[Sal]{_Salamon_} 
S. Salamon, 
Quaternionic K\"ahler manifolds,
Inv. Math. {\bf 67} (1982) 143--171.


\bibitem[Saw]{_Sawon_} 
Sawon, J. 
{\em Abelian fibred holomorphic symplectic
  manifolds}, Turkish
Jour. Math. 27 (2003), no. 1, 197-230, math.AG/0404362.

\bibitem[SYZ]{_SYZ:MS_is_T_du_}  
A. Strominger, S.-T. Yau, and E. Zaslow, 
{\em Mirror Symmetry is T-duality,} Nucl. Phys.
     B479, (1996) 243-259.




\bibitem[V1]{_Verbitsky:coho_announce_} 
Verbitsky, M.,
{\it Cohomology of compact hyperk\"ahler manifolds
and its applications,} alg-geom 
electronic preprint 9511009, 12 pages, LaTeX,
also published in: GAFA vol. 6 (4) pp. 601-612 (1996).


\bibitem[V2]{_Verbitsky:cohomo_} 
Verbitsky, M., {\it Cohomology of 
compact hyperk\"ahler manifolds,}  alg-geom 
electronic preprint 9501001, 89 pages, LaTeX.


\bibitem[V3]{_Verbitsky:SYZ_}
Verbitsky, M., 
{\em Hyperkahler SYZ conjecture and semipositive line bundles},
 arXiv:0811.0639, GAFA 19, No. 5 (2010) 1481-1493.


\bibitem[V4]{_Verbitsky:parabolic_}
Misha Verbitsky, 
{\em Parabolic nef currents on hyperk\"ahler manifolds},
arXiv:0907.4217, 22 pages. 


\bibitem[V5]{_V:Torelli_}
Verbitsky, M.,
{\em A global Torelli theorem for hyperk\"ahler manifolds,}
arXiv: 0908.4121,
 Duke Math. J. Volume 162, Number 15 (2013), 2929-2986.


\bibitem[Y]{_Yau:Calabi-Yau_} 
S. T. Yau,  {\em On the Ricci curvature of a compact K\"ahler manifold 
and the complex Monge-Amp\`ere equation I,}  Comm. on Pure and Appl.
Math. 31, 339-411 (1978).



\end{thebibliography}
\end{document}